\documentclass[16pt]{article}
\title{Understanding Stochastic Differential Equations}
\author{Pat Muldowney}

\newtheorem{definition}{Definition}

\newtheorem{example}{Example}

\newcommand{\R}{\mathbf{R}}

\newcommand{\D}{\mathcal{D}}

\newcommand{\Pa}{\mathcal{P}}

\newcommand{\ve}{\varepsilon}

\usepackage[pdftex]{graphicx}

\begin{document}
\date{}\maketitle
\begin{abstract}
This essay explores the meaning of stochastic differential equations and stochastic integrals.
It sets these subjects in a context of Riemann-Stieltjes integration. It is intended as a comment or supplement to \cite{MTRV}.
\end{abstract}

It\^{o}'s formula is an example of a stochastic differential equation:
\begin{equation}\label{ito1}
dY_s = \frac{\partial f(X_s)}{\partial s}ds + \frac 12 \frac{\partial ^2 f(X_s)}{\partial X_s^2 }ds + \frac{\partial f(X_s)}{\partial X_s} dX_s.
\end{equation}
For $Y_s = f(X_s)$, this formula is an evocative, shorthand way of writing
\begin{equation}\label{ito2}
Y_t - Y_0 =\int_T dY_s =
\int_T \frac{\partial f(X_s)}{\partial s}ds + \frac 12 \int_T  \frac{\partial ^2 f(X_s)}{\partial X_s^2 }ds
+ \int_T \frac{\partial f(X_s)}{\partial X_s} dX_s.
\end{equation}
If the various expressions in this equation represented ordinary numbers and functions,
then the presence in the equation of various integration processes
might incline us to call (\ref{ito2}) an integral equation.

But, while the symbol ``$f$'' in both equations is actually an ordinary, deterministic function (such as the operation of
taking the square of some operand), the symbols $X$ and $Y$ do not represent ``ordinary'' functions or definite numbers.
Instead, they are ``random variables'', that is,  quantities which are indefinite or unknown, to the extent
that they can be predicted only within some margin of error.

The presence of ``$=$'' in the equation indicates that it is an exact statement about actual quantities.
It\^{o}'s formula can be best regarded as an exact statement about margins of error in uncertain quantities.

In other words, it deals with probability distributions of unpredictable quantities which are obtained by means
of various operations in the formula, such as the integration operation. So It\^{o}'s formula
 can be regarded as a kind of integral equation in which the integrals are the type known as \emph{stochastic integrals}.

What is a stochastic integral? What is the meaning of It\^{o}'s formula? These questions are not trivial.
They can be answered in a loose or intuitive manner, but deeper and more exact understanding can be challenging for non-specialists.
And since they are a fundamental part of many important practical subjects, such as finance and communications,
an understanding of them which is merely loose or hazy can be a serious barrier to competent practice in such subjects.

The purpose of this essay is to outline an introduction to stochastic integrals
which is less difficult than the standard textbook treatment of this subject.
It uses Henstock's non-absolute integration instead of Lebesgue integration.
The essay  explores, compares, and contrasts these two methods of integration, with a view to assessing their role in stochastic integrals.

\section{Random variables}
Broadly speaking---at the risk of haziness and looseness!---a \emph{random variable} is a mathematical representation of a
measurement (an experiment,  trial, or observation) of some uncertain or unpredictable occurrence or value.
For instance, the random variable $Z$ could represent a single throw of a die,
so $Z$ represents possible outcomes $\{z=1, \ldots , z=6\}$ with probabilities $\{ \frac 1 6,\ldots ,\frac 16\}$.
Or it could represent measurement of a standard normal variable whose possible values are the real numbers $z \in \R$,
with standard normal probability distribution $\mathbf{N}(0,1)$.

Suppose the throw of the die yields a payoff or outcome $y=f(z)$ obtained by the following deterministic calculation:
\[
y = \left\{
\begin{array}{rl}
-1 & \mbox{if }  z=1 \\
+1 & \mbox{if }  z=6 \\
0 & \mbox{otherwise.}
\end{array}
\right.
\]
This particular experiment or game depends on (is \emph{contingent} on) the outcome of the experiment $Z$,
and can be denoted by $Y=f(Z)$. Where $Z$ has six possible outcomes, with a uniform probability distribution,
$Y$ has three possible outcomes whose probability distribution can easily be deduced by means of the deterministic calculation $f$.
The probability distribution\footnote{The probability distribution (``margin of error'')
carries the essential information specifying the character of the random variable or experiment.
It is often convenient to include other ``potential'' values or outcomes which are \textbf{not} actually possible or ``potential''.
For instance, in the die-throwing experiment we can declare that every real number is a potential outcome. In that case we assign probability
zero to the impossible outcomes. This does not change the random variable or its probability distribution in any essential way
that affects its mathematical meaning} of $Y$ is $y=-1$ with probability $\frac 1 6$, $y=+1$ with probability $\frac 1 6$, and $y=0$
with probability $\frac 2 3$.

We can easily invent such \emph{contingent} random variables or gambling games using more than one throw of the die,
and with payoff $Y$ dependent on some calculation based on the joint outcome  of the successive throws.

This intuitive formulation is compatible with the formal and rigorous conception of a random variable as a $P$-measurable function
whose domain is a $P$-measurable sample space $\Omega$. This twentieth century injection of mathematical rigor brought
about a great extension of the depth and scope of the theory of probability and random variation,
including the development of many new spheres of application of the theory.

These applications often involve \emph{stochastic processes}.
Suppose $T$ is some set of indexing elements $\{s\}$. For instance, $T$ could be an interval of real numbers $[a,b]$.
A stochastic process $Y=Y_T$ is a family $Y=\left(Y(s)\right)_{s \in T}$, for which each element $Y(s), = Y_s$, is a
random variable. A sample path $\left(y(s)\right)_{s \in T}$ of the process $Y=Y_T$ can be thought of as a function
$y: T \mapsto \R$ in which, for each $s$, $y(s)$ (or $y_s$) is a possible outcome of the random variable (measurement, experiment, trial) $Y(s)$.

\section{Stochastic integrals}
Take $T=[0,t]$.
Equation (\ref{ito2}) above appears to be the result of applying an integration operation $\int_T$ to the equation (\ref{ito1}).
If this is the case, and if this step is justified, then comparison of
(\ref{ito1}) and (\ref{ito2}) implies (without delving into their actual meaning) that
\begin{equation}\label{most basic stochastic integral}
\int_T dY_s,=\int_0^t dY_s,=Y_t - Y_0;\;\;\mbox{or }\;\int_T dY(s),=\int_0^t dY(s),=Y(t) - Y(0).
\end{equation}
On the face of it, a clear and precise understanding of this simplest of all possible stochastic integrals would seem to be the \emph{sine
qua non} of this subject. Expressed as a stochastic differential equation, it is the tautology $dY_s = dY_s$.
Whatever (\ref{most basic stochastic integral}) actually means, it seems consistent enough with more familiar forms of integration
of the Stieltjes kind, in the somewhat loose and uncritical sense that the integral (or sum) of increments $dY$ gives an overall increment.

Advancing a little bit further, take a deterministic function $f$, and
consider $\int_T f(Y_s)dY_s$ (or $\int_T f(Y(s))dY(s)$), which is a more general version of        $\int_T dY_s$.
If  $y$ is a sample path of the process $Y$, the expression
\begin{equation}\label{first stochastic integral}
\int_T f(y(s)) dy(s)\;\;\;\mbox{ or }\;\;\int_T f(y_s) dy_s
\end{equation}
is a Stieltjes-type integral, which, if it exists, may be thought of as some limit of Riemann sums
$\sum f(y(s)) \Delta y(s)$ or $\sum f(y(s_j)) \left(y(t_j) - y(t_{j-1})\right)$,
where the finite set of points $t_j$ form a partition of the interval $T=[0,t]$, with
$t_{j-1} \leq s_j \leq t_j$ for each $j$.

From the point of view of basic mathematical analysis, unlike (\ref{most basic stochastic integral}) which is about ``margins of error'' in
probabilistic measurement,
there is nothing problematic about (\ref{first stochastic integral})---this Riemann-Stieltjes-type integral may or may not exist
for particular functions $y$ and $f$, but it is a fairly familiar subject for anyone who has studied basic Riemann-type integration.

In the Riemann sums for (\ref{first stochastic integral}),
{some applications require that $s_j =t_{j-1}$ for each $j$. Cauchy's approach to the theory of integration
used approximating sums with $s_j = t_{j-1}$ or $s_j = t_{j}$, so such sums can be called \emph{Cauchy sums} rather than Riemann sums.}
In any event, there are various ways, including the Lebesgue method,
in which we can seek to define an integral $\int_T f(y(s)) dy(s)$ for sample paths $y_T=(y(s))_{s \in T}$ of a stochastic process $Y = Y_T$.

Suppose a Stieltjes-type integral of $f(y(s))$ is calculated with respect to the increments $y(I) : = y(t_j) - y(t_{j-1})$ of the function $y_T$.
For instance, if $f$ is a function taking some fixed, real, constant value such as $1$, then a ``naive'' Riemann sum calculation
on the domain $T=[0,t]$, with $t_0=0$ and $t_n=1$
gives
$
\sum f(y(s)) y(I) \;\;=\;\;\sum_{j=1}^n y(I) \;\;=$
\[
 = (\left(y(t_1) - y(0)\right) + (\left(y(t_2) - y(t_1)\right) + \cdots + (\left(y(1) - y(t_{n-1})\right)
 =y(1) - y(0)\]
for \textbf{every} sample outcome $y_T$ of the process $Y_T$.
So it is reasonable---in some ``naive'' way---to claim that, for this particular function $f$,
the Riemann-Stieltjes integral exists for all outcomes $y_T$:
\[
\int_T f(y_s) dy_s =\int_0^t dy(s) = y(t)-y(0).
\]
One might then be tempted\footnote{A warning about this temptation is provided in Example \ref{Dirichlet 1} below.}
to apply such an argument to step functions $f$, and perhaps to try to extend it to some class of continuous functions $f$,
especially if we are only concerned with sample paths $y_T$ which are continuous.

But the key point here is that, given a stochastic process $Y=Y_T$, and given certain
deterministic functions $f$,
real values $\int_T f(y(s)) dy(s) $ can be obtained for each sample path $y=y_T$ by means of a recognizable Stieltjes integration procedure.

Can this class of real numbers or outcomes  be related somehow to some identifiable random variable $Z$
which possesses some identifiable probability distribution (or ``margin of error'' estimates)?

If so, then $Z$ might reasonably be considered to be the random variable obtained by integrating,
in some Stieltjes fashion, the random variable $f(Y_s)$ with respect to the increments $Y(I) = Y(t_j) - Y(t_{j-1})$ of the stochastic process $Y_T$.

In other words, $Z$ is the stochastic integral
$\int_T f(Y_s) dY_s$.

To justify the latter step, a probability distribution (or ``margin of error'' data) 
for $Z$ must be determined. But, in the case of the constant function $f$
given above ($f(y_s)=1$), this is straightforward. Because, with $f(y_s)=1$ for all outcomes $y_s$ in all sample paths (or joint outcomes) $y_T$,
the distribution function obtained for the Riemann sum values
$\sum f(y_s) y(I)$ is simply the known distribution function of the outcomes $y(t) - y(0)$ of the random variable $Y(t) - Y(0)$.

This distribution is the same for all partitions of $T=[0,t]$.
So it is reasonable to take it to be the distribution function of the stochastic integral $Z=\int_T f(Y_s) dY_s$.
For constant $f$ this seems to provide meaning and rationale for (\ref{first stochastic integral}).

What this amounts to is a naive or intuitive interpretation of stochastic integration which seems to hold for some elementary functions $f$.
This approach can be pursued further to give a straightforward interpretation---indeed, a ``proof''---of It\^{o}'s formula,
at least for the unchallenging functions $f$ mentioned above.

But what of the standard or rigorous theory of stochastic integration?

\section{Standard theory of stochastic integration}
Unfortunately, this theory
 cannot accommodate the naive or intuitive construction of the simple stochastic integrals described in the preceding section.
 Broadly speaking, the elementary Riemann sum type of calculation is not adequate for the kinds of analysis needed in this subject.
 It is not possible, for instance, to apply a monotone convergence theorem, or a dominated convergence theorem,
 to simple Riemann and Riemann-Stieltjes integrals. Historically, these kinds of analysis and proof have been supplied by Lebesgue-type integrals
 which, while requiring a measure function as integrator, cannot be simply defined by means of the usual arrangement\footnote{But Section
 \ref{-Complete integration}  shows that Lebesgue integrals are essentially Riemann-Stieltjes integrals.}
 of Riemann sums.

And this is where the difficulty is located. Suppose, for instance, that the stochastic process $Y_T$ that we are dealing with is a
standard Brownian motion. In that case any sample path $y_T$ is, on the one hand, almost surely continuous---which is ``nice'';
but, on the other hand, it is almost surely
\emph{not} of bounded variation in every interval $J$ of the domain $T=[0,1]$. And the latter is ``nasty''.

This turns out to be very troublesome if we wish to construct a Lebesgue-Stieltjes integral
using the increments $y(I), = y(t_j)-y(t_{j-1})$, of a sample path which is continuous but not of bounded variation in any interval.

The problem is that, in order to construct a Lebesgue-Stieltjes measure from the increments $y(I)$, we must separate the
non-negative increments $y_+(I)$ from the negative-valued increments $y_-(I)$,
\[
y(I) = y_+(I) - |y_-(I)|,
\]
and try to construct a non-negative measure from each of the components. But, because $y$ is not of bounded variation,
the construction for each component diverges to infinity on every interval $J$. Thus the standard theory of stochastic
integration encounters a significant difficulty at the very first step (\ref{first stochastic integral}).

To summarize:
\begin{itemize}
\item In the standard It\^{o} or Lebesgue integral approach, the most basic calculation of the integral of a constant function $f(Y_T)$,
with respect to the increments $dY$ of a Brownian process, fails because the Lebesgue-Stieltjes measure does not exist.
\item
On the other hand, if Riemann sums of the increments of the process $Y_T$ are used, then, by cancelation, a finite result is obtained
for each Riemann sum---a result which agrees with what is intuitively expected.
\end{itemize}
In the standard Lebesgue (or It\^{o}) theory of stochastic integration---in \cite{O} for instance---this
problem is evaded by \emph{postulating} a finite measure
$\mu_y(J)$ for each sample path, and then constructing a weak form of integral which, in the case of Brownian motion,
is  based on certain helpful  properties of this process.

The trouble with this approach is that it produces a quite  difficult theory
which does not lend itself to the natural, intuitive interpretation described above.

However, elementary Riemann-sum-based integration is not generally considered to have the analytical power possessed by Lebesgue-style integration.
And a great deal of analytical power is required in the theory of stochastic processes. So at first sight it seems that
we are stuck with the standard theory of stochastic integration, along with all its baggage of subtlety and complication.

But this is not really the case. The good news is that is actually possible to formulate the theory of stochastic integrals
using Riemann sums instead of the measures of Lebesgue theory.

\section{Integration of functions}
To see this, it is first necessary to review the various kinds of integration which are available to us.

First consider the basic Riemann integral, $ \int_a^b f(s) ds$, of a real-valued, bounded, continuous function $f(s)$ on an interval $[a,b]$.
Let $\Pa$ be a partition of $[a,b]$;
\[
\Pa: \;\;\;a=t_0<t_1<t_2< \cdots < t_n=b,
\]
for any choice of positive integer $n$ and any choice of $t_j$, $1\leq j <n$.
For any $u<v$ and any interval $I$ with end-points $u$ and $v$, write $|I| =v-u$. Denoting intervals $]t_{j-1}, t_j]$ by $I_j$ let
\[
U_\Pa=\sum_{j=1}^n  P_j|I_j|,\;\;\;\;\;\;L_\Pa=\sum_{j=1}^n p_j|I_j|
\]
where
\[
P_j = \sup \{f(s): s \in I_j \},\;\;\;\;\;\; p_j = \inf \{f(s): s \in I_j \}.
\]
\begin{definition}\label{def riemann}
Define the \emph{upper Riemann integral} of $f$ by
\[
U:=\inf\{L_\Pa: \mbox{ all partitions } \Pa \mbox{ of }[a,b]\},
\]
and the \emph{lower Riemann integral} of $f$ by
\[
L:=\sup\{l_\Pa: \mbox{ all partitions } \Pa \mbox{ of }[a,b]\}.
\]
Then $U_\Pa \geq L_\Pa$ for all $\Pa$, and if $U= L$ we say that $f$ is \emph{Riemann integrable}, with
\[
\int_a^b f(s) ds :=U=L.
\]
\end{definition}
Write the partition $\Pa $ as $ \{I\}$ where each $I$ has the form $I_j =\, ]t_{j-1}, t_j]$, with $|I_j| = t_j-t_{j-1}$, and Riemann sum
\[
(\Pa)\sum f(s)|I| = \sum_{j=1}^n f(s_j)|I_j|.
\]
Suppose $g(s)$ is a real-valued, monotone increasing function of $s \in [a,b]$, so $g(s) \geq g(s')$ for $s>s'$.
For any interval $I$  with end-points $u$ and $v$ ($u<v$), define the increment or interval function $g(I)$ to be $g(v)-g(u)$.
\begin{definition}\label{def riemann-stieltjes 1}
If $|I|$ and $|I_j|$ are replaced by $g(I)$ and $g(I_j)$ in Definition \ref{def riemann} of the Riemann integral, then the resulting integral is called
the \emph{Riemann-Stieltjes integral} of $f$ with respect to $g$, $\int_a^b f\;dg$ or $\int_a^b f(s) dg(s)$.
\end{definition}
In fact if we start with the latter definition the Riemann integral is a special case of it,
obtained by taking the point function $g$ to be the identity function $g(s) =s$.

If $g(s)$ has bounded variation it can be expressed as the difference of two monotone increasing, non-negative point functions,
\[
g(s) = g_+(s) - (-g_-(s)),
\]
and the Riemann-Stieltjes integral of $f$ with respect to $g$ can then be defined as the
difference of the Riemann-Stieltjes integrals of $f$ with respect to $g_+$ and $-g_-$, respectively.

The following result is well known:
if real-valued, bounded $f$ is continuous and if real-valued $g$ has bounded variation then $\int_a^b f\;dg$ exists.

As suggested earlier, the  Lebesgue integral of a real-valued point function $k$ with respect to a measure $\mu$ can be viewed,
essentially, as a Riemann-Stieltjes integral in which the
point-integrand  $k(\omega)$ satisfies the condition of \emph{measurability}.
To explain this statement further,
consider a measure space $(\Omega, \mathcal{A}, \mu)$ with non-negative measure $\mu$ on a sigma-algebra $\mathcal{A}$
of $\mu$-measurable subsets of the arbitrary measurable space $\Omega$. Thus, if $\mu(\Omega) =1$, the measure space is a probability space.
Suppose the point-integrand $k$ is a bounded real-valued $\mu$-measurable function on the domain $\Omega$.
Then there exist real numbers $c$ and $d$ for which
\[
c \leq k(\omega) \leq d \;\;\;\mbox{ for all }\;\;\;\omega \in \Omega.
\]
Also, for each sub-interval $J$ of $[c,d]$, measurability of $k$ implies $\mu(k^{-1}(J))$ is defined.
The basic definition of the Lebesgue integral of $k$ with respect to $\mu$ on $\Omega$ is as follows.
\begin{definition}\label{def lebesgue}
Let $\mathcal{Q} = \{J_j\} =\{]v_{j-1},v_j]\}$
be a partition of $[c,d]$,
\[
\mathcal{Q}: \;\;\;c=v_0<v_1<v_2 < \cdots < v_n=d,
\]
and let
\[
L_\mathcal{Q}= \sum_{j=1}^n v_{j-1}\mu(k^{-1}(J_j)),\;\;\;\;\;\;U_\mathcal{Q}= \sum_{j=1}^n v_{j}\mu(k^{-1}(J_j)).
\]
Let $L:= \sup\{L_\mathcal{Q}: \mathcal{Q}\}$, $U:= \inf\{U_\mathcal{Q}: \mathcal{Q}\}$, the supremum and infimum being taken over all
partitions $\mathcal{Q}$ of $[c,d]$.
If $L=U$, then their common value is the \emph{Lebesgue integral} $\int_\Omega k(\omega)d\mu$.
\end{definition}
An advantage of Lebesgue integration over Riemann integration is that the former has theorems,
such as the dominated and monotone convergence theorems which,
under certain condition, make it possible for instance to change the order of integration and differentiation.
Also, Fubini's and Tonelli's theorems allow exchange of order of multiple integrals.

What makes  ``good'' properties such as these possible is \emph{measurability} of the integrand $k$.
But the Lebesgue integral itself is, by definition, a Riemann-Stieltjes-type integral.
To see this, for each $u \in [c,d]$ define the monotone increasing function
\begin{equation}\label{lebesgue as riemann-stieltjes (1)}
g(u) = \mu\left(k^{-1}([c,u])\right),
\end{equation}
and take the point function $h(u)$ to be the identity function $h(u)=u$. Then the construction\footnote{The integral of a point function $h(u)$ with respect to a point function $g(u)$ can be addressed either as a Riemann-Stieltjes construction or as a Lebesgue-Stieltjes construction.  When $h(u)=u$ and $g(u) = \mu\left(k^{-1}([c,u])\right)$ the former approach gives the Lebesgue integral $\int_\Omega k(\omega) d\mu$. On the other hand, if the Lebesgue-Stieltjes construction is attempted with $h(u)=u$ and $g(u) = \mu\left(k^{-1}([c,u])\right)$, we simply replicate the Riemann-Stieltjes construction of the Lebesgue integral $\int_\Omega k(\omega) d\mu$, and nothing new emerges.}
in Definition \ref{def lebesgue}
 shows that
\begin{equation}\label{lebesgue as riemann-stieltjes (2)}
\int_\Omega k(\omega) d\mu \;\;\;=\;\;\;\int_c^d h(u) \;dg(u),\;\;\;=\;\;\;\int_c^d u\;dg.
\end{equation}
In other words, when combined with the measurability property of the point-integrand,
this particular Riemann-Stieltjes construction gives the ``good'' properties required in the integration of functions.

\section{Riemann definition}\label{Riemann definition}
But in fact a Riemann construction can give these ``good'' properties \textbf{without} postulating measurability in the
definition\footnote{And if measurability is
redundant in the definition, then so is the  measure space structure.}  of the integral.
To see this, we start again by considering a more general and more flexible definition of basic Riemann and Riemann-Stieltjes integration
which generalizes the  construction of these integrals as given above in Definitions \ref{def riemann} and \ref{def riemann-stieltjes 1}.

The proposed, more general, definition of the Riemann-Stieltjes integral is applicable to real- or complex-valued functions $f$ (bounded or not);
and to real- or complex-valued functions $g$, with or without bounded variation.
\begin{definition}\label{def riemann-stieltjes 2}
The function $f$ is Rie\-m\-ann-Stieltjes integrable with respect to $g$, with integral $\alpha$, if, given $\ve>0$,
there exists a constant $\delta>0$ such that, for every partition $\Pa = \{I\}$ of $[a,b]$ satisfying $|I| < \delta$ for each $I \in \Pa$,
the corresponding Riemann sum satisfies
\[
\left| \alpha - (\Pa)\sum f(s)g(I) \right| < \ve,
\]
so $\alpha = \int_a^b f\;dg$.
\end{definition}
If $g$ is the identity function $g(s)=s$ then Definition \ref{def riemann-stieltjes 2} reduces to the ordinary Riemann integral of $f$, $\int_a^b f(s)ds$.

{Definition \ref{def riemann-stieltjes 2} does not embody conditions which
ensure the existence of the integral. Such integrability conditions are not postulated but are deduced, in the form of theorems,
from the definition of the integral.}

Thus, if the function properties specified, respectively, in Definitions \ref{def riemann},
\ref{def riemann-stieltjes 1}, and \ref{def lebesgue} above are assumed, the integrability in each case follows from
Definition \ref{def riemann-stieltjes 2}; and Definitions \ref{def riemann},
\ref{def riemann-stieltjes 1}, and \ref{def lebesgue} become theorems of Riemann, Riemann-Stieltjes, and Lebesgue integration, respectively.

Definition \ref{def lebesgue} can now be expressed in terms of Definition \ref{def riemann-stieltjes 2}, using the formulations
(\ref{lebesgue as riemann-stieltjes (1)}) and
(\ref{lebesgue as riemann-stieltjes (2)}), and assuming measurability of the integrand $f$ with respect to measure space $(\Omega, \mathcal{A}, \mu)$.
\begin{definition}\label{def lebesgue 2}
The function $f$ is Lebesgue integrable with respect to measure $\mu$, with integral $\int_\Omega f(\omega)d\mu=\alpha$, if, given $\ve>0$,
there exists a constant $\delta>0$ such that, for every partition $\mathcal{Q} = \{J\}$ of $[c,d]$ satisfying $|J| < \delta$ for each $J \in \mathcal{Q}$,
the corresponding Riemann sum satisfies
\[
\left| \alpha - (\mathcal{Q})\sum h(u)g(J) \right| < \ve,
\]
where $h(u)=u$ is the identity function on $[c,d]$;
so $\alpha = \int_c^d h(u)dg(u)= \int_c^d u\;dg$.
\end{definition}
Thus, by definition, the Lebesgue integral $\int_\Omega f(\omega)d\mu$, with domain $\Omega$, is the Riemann-Stieltjes integral $\int_c^d u\;dg$,
with domain $[c,d]$.

The following result is an obvious
consequence of {Definition} \ref{def riemann-stieltjes 2}.
 If $f$ has constant value $\beta$ and if $g$ is an arbitrary real- or complex-valued function,
then $\int_a^b f\;dg$ exists and equals $\beta(g(b) - g(a))$.
This follows directly from Definition \ref{def riemann-stieltjes 2} since, for every partition $\Pa$ of $[a,b]$, cancelation of terms gives
\[
(\Pa)\sum f(s)g(I) = \beta\sum_{j=1}^n g(t_j) - g(t_{j-1}) = \beta\left(g(b) - g(a)\right).
\]
This result does not in general hold for Lebesgue-Stieltjes integration, as the latter
requires that $g(s)$ be resolved into its negative and non-negative components, $g(s) = g_+(s) - (-g_-(s)$,
and convergence may fail when the integral is calculated with respect to each of these components separately.

Example \ref{Dirichlet 1} below shows that, though constant functions $f$ are Riemann-Stieltjes integrable with respect to any integrator function $g$,
this does not necessarily extend to step functions $f$.

Definition \ref{def riemann-stieltjes 2} of the Riemann or Riemann-Stieltjes integral
does not postulate any boundedness, continuity, measurability
or other conditions for the integrand $f$. But, as already stated, in the absence of integrand measurability
and the construction in Definition \ref{def lebesgue}, this method of integration does not deliver
 good versions of monotone and dominated convergence theorems, or Fubini's theorem.

\section{-Complete integration}\label{-Complete integration}
Developments in the subject since the 1950's---developments which were originated independently by R.~Henstock
and J.~Kurzweil---have made good this deficit in the basic Riemann and Riemann-Stieltjes construction. In this new development of the subject,
Definition \ref{def riemann-stieltjes 2} of the Riemann-Stieltjes integral is amended
as follows.

\begin{definition} \label{def stieltjes-complete}
A function $f$ is \emph{Stieltjes-complete} integrable with respect to a function $g$, with integral $\alpha$ if, given $\ve>0$,
there exists a function $\delta(s) >0$ such that
\[
\left| \alpha - (\Pa)\sum f(s)g(I) \right| < \ve
\]
for every partition $\Pa$ such that, in each term $f(s)g(I)$ of the Riemann sum, we have $s-\delta(s) <u \leq s \leq v < s+\delta(s)$,
where $u$ and $v$ are the end-points of the partitioning interval $I$.
\end{definition}
In other words, where $|I|$ is less than a constant $\delta$ in the basic Riemann-Stieltjes definition,
we have $|I| < \delta (s)$ in the new definition.
Write $\alpha = \int_{[a,b]} f(s)g(I)$, or $\int_{[a,b]} f\;dg$, for the Stieltjes-complete integral whenever it exists.

Again, if the integrator function $g$ is the identity function $g(s)=s$, the resulting integral (corresponding to the basic Riemann integral), is the
\emph{Riemann-complete} integral of $f$, written $\alpha = \int_{[a,b]}f(s)|I|$, or $\int_{[a,b]}f(s)ds$.
The latter is also known as the Henstock integral,
the Kurzweil integral, the Henstock-Kurzweil, the generalized Riemann integral, or the \emph{gauge} integral since in this context the function
$\delta(s)>0$ is called a gauge.

It is obvious that every Riemann (Riemann-Stieltjes) integrable integrand is also Riemann-complete (Stieltjes-complete) integrable,
as the gauge function $\delta(s)>0$ of Definition \ref{def riemann-stieltjes 2} can be taken to be the constant $\delta >0$ of
Definition \ref{def riemann} and Definition \ref{def riemann-stieltjes 1}.

This argument  indicates a \emph{Lebesgue-complete} extension of the Lebesgue integral, by replacing the
constant $\delta >0$ of Definition \ref{def lebesgue 2} with a variable gauge $\delta(u)>0$:
\begin{definition}\label{def lebesgue-complete}
 Let $h(u)=u$ be the identity function on $[c,d]$.
The function $f$ is \textbf{Lebesgue-complete} integrable with respect to measure $\mu$, with integral $\int_\Omega f\;d\mu=\alpha$, if, given $\ve>0$,
there exists a gauge $\delta(u)>0$ for $c \leq u \leq d$, such that
\[
\left| \alpha - (\mathcal{Q})\sum h(u)g(J) \right| < \ve,
\]
for every partition $\mathcal{Q} = \{J\}$ of $[c,d]$ satisfying
\[u-\delta(u)<v_{j-1} \leq u \leq v_j < u+\delta(u)\]
for each $J =\,]v_{j-1},v_j]\in \mathcal{Q}$.
\end{definition}
In that case $\alpha = \int_{[c,d]} h(u)g(J)= \int_{[c,d]} u\;g(J)$, and the
 Lebesgue-complete integral is a special case of the Stieltjes-complete integral---a special case in which a measure space structure exists and
for which the integrand is measurable.
So it is again clear  that every Lebesgue integrable integrand is Lebesgue-complete integrable; since the former is, in effect,
a Riemann-Stieltjes integral, the latter is a Stieltjes-complete integral,
and every Riemann-Stieltjes integrable function is also Stieltjes-complete integrable. (No special notation has been introduced here
to distinguish the Lebesgue integral $\int_\Omega f\;d\mu$ from its Lebesgue-complete counterpart.)

If the measurable domain $\Omega$ is a real interval such as $[a,b]$, then some ambiguity arises in the interpretation of the Lebesgue integral
as an integral of the gauge, or generalized Riemann, kind. The reason for the ambiguity is as follows. Assuming the existence of the
Lebesgue integral $\int_{\Omega}f(\omega) d\mu, =\int_{[a,b]}f(\omega) d\mu$,
where $\omega$ now represents real numbers in the domain $[a,b]$, then we are assured of
the existence of the Stieltjes and Stieltjes-complete (or Lebesgue-complete) integrals $\int_c^d u\;dg$ and $\int_{[c,d]}u\;g(J)$,
respectively, with
\[
\int_{[a,b]}f(\omega) d\mu = \int_c^d u\;dg = \int_{[c,d]}u\;g(J),
\]
where the values $u, = h(u)$, are elements of $[c,d]$ and $h$ is the identity function on $[c,d]$.

But in this case, letting $\omega = s$ denote points of the domain $[a,b]$ and with $I$ denoting subintervals of $[a,b]$,
the function $\mu(I)$ is defined on intervals $I$, and two different Stieltjes-type constructions are possible.

First, there is the Riemann-Stieltjes integral $\int_c^du \;dg$ which defines the Leb\-esgue integral $\int_{\Omega}f(\omega) d\mu,
= \int_{[a,b]}f(\omega) d\mu$. Secondly, there is the gauge integral $\int_{[a,b]}f(s) \mu(I)$ which has a Stieltjes-complete construction.

It is then meaningful to consider whether, with $f$ measurable, existence of the Lebesgue integral
$\int_{[a,b]}f(\omega) d\mu$ implies existence of the Stieltjes-complete integral $\int_{[a,b]} f(s) \mu (I)$, and whether
\[
\int_c^du \;dg = \int_{[a,b]}f(s) \mu(I)
\]
holds,\footnote{There is a considerable literature on this question, which is usually answered as:
\emph{``Every Lebesgue integrable function on an interval of the real numbers $\R$ is also
Henstock-Kurzweil integrable.''} If the domain of the integrand is a measurable space $\Omega$ which is \textbf{not}
a subset of $\R$ or $\R^n$, then the appropriate way to formulate the corresponding Henstock-Kurzweil (or -complete)
integral is in the form $\int_{[c,d]} u\,g(J)$
described in Definition \ref{def lebesgue-complete}.}
 the first of these integrals being the  Lebesgue integral
$\int_{[a,b]}f(\omega) d\mu $, which, by Definition \ref{def lebesgue 2}, is interpreted as the Riemann-Stieltjes integral $\int_c^du \;dg$.

To see that these two integrals coincide, take $f$ to be a bounded, measurable function on $[a,b]$. This can be expressed as the difference of two
non-negative, bounded, measurable functions $f_+$ and $f_-$. Accordingly, and without loss of generality,
take $f$ to be non-negative, bounded, measurable. Then the Lebesgue integrable
function $f$ is the $\mu$-almost everywhere point-wise limit of a monotone increasing sequence of step functions $f_j$.
With $\omega = s$,
each step function $f_j$ is
Lebesgue integrable, with Lebesgue integral $\int_{[a,b]}f_j(\omega) \;d\mu$; and
each step function $f_j$ is Stieltjes-complete integrable, with Stieltjes-complete integral
$\int_{[a,b]}f_j(s) \mu(I)$, and
\[
\int_{[a,b]}f_j(\omega) \;d\mu = \int_{[a,b]}f_j(s) \mu(I)
\]
for each $j$. (This statement is also true if ``Lebesgue integral'' and ``Lebesgue integrability'' are replaced by ``Lebesgue-complete integral'' and ``Lebesgue-complete integrability''.)

By the monotone convergence theorem of Lebesgue integration (or, respectively, by the monotone convergence theorem of Lebesgue-complete integration),
\[
\int_{[a,b]}f_j(\omega) d\mu \rightarrow \int_{[a,b]}f(\omega)  d\mu
\]
as $j \rightarrow \infty$. By the monotone convergence theorem of Stieltjes-complete integration, $f(s)\mu(I)$ is Stieltjes-complete integrable and
\[
\int_{[a,b]}f_j(s) \mu(I) \rightarrow
\int_{[a,b]}f(s) \mu(I)
\]
as $j \rightarrow \infty$. Since corresponding integrals of the pair of sequences are equal, their limits are equal:
\[
\int_{[a,b]}f(\omega)  d\mu = \int_{[a,b]}f(s)\mu(I).
\]
This is the gist of a proof that existence of a Lebesgue integral (or of a Lebesgue integral) on a real domain implies existence of the corresponding Stieltjes-complete integral
on the same domain,
and equality of the two.

Thus the above argument can be applied to either the Lebesgue or the Lebesgue-complete integral on $\Omega = [a,b]$ in conjunction, respectively, with the corresponding Stieltjes-complete integral
on the same domain. In effect, if the domain $\Omega$ is a subset of $\R$, and if $f$ is Lebesgue integrable or Lebesgue-complete integrable with respect to $\mu$,
then $f(s)\mu(I)$ is also Stieltjes-complete integrable and the two integrals are equal.

The specific properties of the Lebesgue-complete integral have not been investigated.

As mentioned earlier, constant functions $f$ are Riemann-Stieltjes integrable,
and hence Stieltjes-complete integrable, with respect to any integrator function $g$.
But as the following counter-example shows, this does not necessarily extend to step functions $f$, or any other functions which are not constant.

\begin{example}
\label{Dirichlet 1}
\textbf{Dirichlet function:} For $0 \leq s \leq 1$ let $D(s)$ be $1$ if $s$ is rational, and $0$ otherwise. For $I=\,]u,v]$ let
$D(I) = D(v) - D(u)$. Let $D([0,v])= D(v) - D(0)$. The point function $D(s)$ is discontinuous everywhere, and has infinite variation on every interval
$J \subseteq [0,1]$. If $f(s)$ is constant for $0\leq s \leq 1$, then the Riemann-Stieltjes integral $\int_0^1 f(s)\,dD$ exists and equals
$D(1)-D(0)$; that is, $\int_0^1 f(s)\,dD=0$. \textbf{But if $f$ is not constant on $[0,1]$, then the Riemann-Stieltjes integral of $f$
with respect to $D$
does not exist.}
What about Stieltjes-complete integrability of $f(s)D(I)$?
In fact, if $f$ is not constant on $[0,1]$, then the Stieltjes-complete integral of $f$
with respect to $D$
does not exist.
This is proved in Theorem 1 of \cite{Henstock 1955}, and the proof is reproduced in Theorem 67 of \cite{MTRV}.
\textbf{Thus $f(s) D(I)$ is Riemann-Stieltjes integrable and
Stieltjes-complete integrable on $[0,1]$ if and only if $f(s)$ is constant for $0 \leq s \leq 1$.}
\end{example}
 Historically this is  the first published result (Theorem 1 of \cite{Henstock 1955}) in the theory of -com\-plete integration.

\section{-Complete approach to stochastic integrals}
Returning to stochastic integrals, the -complete method of integration allows us to construct Stieltjes-type Riemann sums for highly oscillatory
expressions which include both positive
and negative terms. Cancelation of terms can occur in the Riemann sum approximations, so the possibility of  convergence is
preserved by this construction.

The Lebesgue construction, on the other hand, requires integral convergence, separately and independently, of the positive and negative components
of the integrand. The difficulty this presents is illustrated in the alternating or oscillating series $\sum_{j=1}^\infty (-1)^{-j} j^{-1}$. If the positive
and negative terms of the series are considered as two separate series then each of them diverges. But the series itself is conditionally
(or non-absolutely) convergent. Similarly, for sample paths $y(s)$ of a stochastic process $Y_T$ the integral $\int_{[0,t]}dy(s)$
does not generally exist when considered as a Lebesgue-Stieltjes integral. But it exists for all sample paths $y_T$, with value $y(t) - y(0)$, when
considered as a Stieltjes-complete integral.

There is no analytical cost or disadvantage in relinquishing the Lebesgue construction in favor of the -complete method. This is because
the important theorems of Lebesgue integration, such as monotone and dominated convergence, are also valid for the -complete approach.
Furthermore, there are other convergence theorems of a similar  kind, specifically designed to deal with highly oscillatory functions
such as those which occur in the theory of stochastic processes but which are beyond the scope of the Lebesgue method.
See \cite{MTRV} for details of these.

However, stochastic integration includes novelties and challenges which have not yet been addressed in this essay.

For Brownian motion processes $X_T$, one of the most important stochastic integrals is $\int_0^t dX_s^2 =t$.
The corresponding integral for a sample path $x(s)$ ($0\leq s \leq t$) is ``$\int_0^t (dx(s))^2$''.
But this expression does not have the familiar form of a Stieltjes-type integral: $\int_a^b f(s)dg$,
which, when $g$ is the identity function, reduces to the even more familiar $\int_a^b f(s)ds$.

In Riemann sum approximation we are dealing with expressions $\sum (x(I))^2$, where, for $I=\,]u,v]$, $x(I) = x(v)-x(u)$.
But traditionally, while a Riemann sum for a Stieltjes integral involves terms $f(s)x(I)$ with integrator function $x(I)$
(in which $f(s)$ can be identically $1$),
we do not usually expect to see integrators such as $(x(I))^2$ or $dX_s^2$.

Another important stochastic integral Brownian motion theory is
\[
\int_0^t X_s dX_s = \frac 12 X_t^2 - \frac 12 t.
\]
For a sample path $x(s)$ of Brownian motion, this involves $\int_0^t x(s) dx(s)$, or, in Riemann sum terms, $\sum x(s)x(I)$.
The latter, as it stands, is a finite sum of terms $x(s)(x(v) - x(u))$ where $I=\,]u,v]$ and $u \leq s \leq v$.
And if we are using the Stieltjes-complete approach as described above, then we might suppose that each $s$ in the Riemann sum
is the special point used in partitions which are constrained by a gauge $\delta(s)$,
\[
s-\delta(s) < u \leq s \leq v < s+\delta(s).
\]
But in fact this is not what is required in the stochastic integral $\int_0^t X_s dX_s$. In Riemann sum format, what is required is
\[
\sum x(u) x(I),\;\;\;\mbox{ or }\;\;\;
\sum x(u) \left(x(v) - x(u)\right),
\]
where the first factor $x(u)$ in the integrand is a point function evaluated at the left hand end-point $u$ of the interval $I=\,]u,v]$.

Sometimes the form $\sum x(w)(x(v)-x(u))$ is used, with $w=u + \frac 12 (v-u)$.

In a way, integrands of form $x(I)^2$, $x(u)x(I)$, or $x(w)x(I)$, are an unexpected innovation.
Their value is calculated from the numbers $u$ and $v$ which specify the interval $I$.
So they can be thought of as functions $h(I)$ of intervals $I$.

But these functions are \textbf{not additive} on intervals.\footnote{If $h(I)$ were finitely additive on intervals $I$ it could
be used to define a point function $h(s):=h([0,s])$, and vice versa. Integrals with respect to finitely additive integrators are
therefore representable as Stieltjes-type integrals, and vice versa.}
In that regard they are unlike the integrators $|I|$ and $x(I)$ which are themselves functions of $I$ but are
finitely additive on intervals, in the sense that, if $J=I_1 \cup \ldots \cup I_n$ is an interval, then
\[
|J| = \sum_{j=1}^n |I_j|,\;\;\;\;\;\;x(J) = \sum_{j=1}^n x(I_j).
\]
Broadly speaking, integration is a summation process in which the summed terms involve functions of intervals.
Up to this point in this essay, the only integrands to be considered included a factor which was an additive function of intervals $I$, such as
the length function $|I|$
or the Stieltjes-type functions $g(I)$ or $x(I)$. But there is nothing inherent in the definition of -complete integrals that requires any
$I$-dependent factor in the integrand
to be additive.

With this in mind, consider again the definition of the -complete integral on an interval $[a,b]$.

Firstly, a \emph{gauge} is a function $\delta(s)>0$, $a\leq s\leq b$. Given $s$, an interval $I=\,]u,v]$
for which $s$ is either an end-point or an interior point, is
 $\delta(s)$\emph{-fine}
if $s-u<\delta(s)$ and $v-s<\delta(s)$. A finite collection $\D = \{(s_1,I_1), \ldots , (s_n,I_n)$ is a \emph{division} of $[a,b]$
if each $s_j$ is either an interior point or end-point of $I_j$ and the intervals $I_j$ form a partition of $[a,b]$.
Given a gauge $\delta$, a division $\D$ is $\delta$-fine if each $(s_j,I_j) \in \D$ is $\delta$-fine.

Now suppose $h$ is a function of elements $(s,I)$. Examples include:
\[
h(s,I) =h_1(I)= |I|,\;\;\;h(s,I) =h_2(s)=s,\;\;\;h_3(s,I)=s^2|I|,\;\;\;h_4(I) = |I|^2.
\]
Given a division $\D=\{(s,I)\}$ of $[a,b]$ whose intervals $I$ form a partition $\Pa$, the corresponding Riemann sum is
\[
(\D)\sum h(s,I), = \sum \{h(s,I): I \in \Pa \}.
\]
\begin{definition}
\label{def Burkill integral}
A function $h(s,I)$ is integrable on $[a,b]$, with integral $\int_{[a,b]}h(s,I) = \alpha$, if, given $\ve>0$.
there exists a gauge $\delta(s)>0$ so that, for each $\delta$-fine division $\D$ of $[a,b]$,
\[
\left| \alpha - (\D)\sum h(s,I) \right| < \ve.
\]
\end{definition}
Applying this definition to the examples, $h_1$ is integrable with integral $b-a$, $h_2$ is not integrable, $h_3$ is integrable with integral
$\frac 13(b^3 - a^3)$, and $h_4$ is integrable with integral $0$. If $h(s,I) = h_5(I) = u^2|I|$ where, for each $I$, $u$ is the left hand end-point of $I$,
then it is not too hard to show that $h_5$ is integrable with integral $\frac 13(b^3 - a^3)$.

Actually, it is the traditional
custom and practice in this subject to only consider integrands $h(s,I) = f(s)p(I)$
where the integrator function $p(I)$ is a measure function or, at least, finitely additive on intervals $I$;
and where the evaluation point $s$
of the point function integrand $f(s)$ is the point $s$ of $(s,I)$  for each $(s,I)\in \D$.
When $p(I) = |I|$, this convention is needed in order to prove the Fundamental Theorem of Calculus.\footnote{The Fundamental Theorem of Calculus
states that if $F'(s)=f(s)$ then $f(s)$ is integrable on $[a,b]$ with definite
integral equal to $F(b)-F(a)$}.

But, while the Fundamental Theorem of Calculus is important in subjects such as differential equations, it hardly figures at all in some other branches
of mathematics such as probability theory or stochastic processes. And we have seen that stochastic integration often requires point integrands
$f(s)$ to be evaluated, not at the points $s$ of $(s,I) \in \D$, but at the left hand end-points of the partitioning intervals $I$.

So, with $I=\,]u,v]$, $f(u)$ is, in fact, an integrand function which depends, not on points $s$ but on intervals $]u,v]$.

These are a few of the ``unexpected innovations'' to be encountered in stochastic integration, giving it a somewhat alien and counter-intuitive feel
to anyone versed in the traditional methods of calculus.

For instance, the stochastic integral $\int_0^t X dX$ is given the value $\frac 12 X(t)^2 - \frac 12 t$ when the process $X(s)$ (with $X(0)=0$)
is a Brownian motion. Introductory treatments of this problem sometimes contrast the expression $\int_0^t X dX$
with the elementary calculus integral $\int x dx$ whose indefinite integral is $\frac 12 x^2$, in which the use of symbols $X$ and $x$
can, in the mind of an inexperienced reader, set up an inappropriate and misleading analogy.

In terms of sample paths, the stochastic integral $\int_0^t X(s) dX(s)$ has representative sample form $\int_0^t x(s) dx(s)$
which is a Stieltjes-type integral with integrator function $x(I), = x(v)-x(u)$, formed from a typically ``zig-zag'' Brownian path $x(s)$,
$0<s\leq t$, with $x(0)=0$. Then the notation for the contrasting elementary calculus integral is not $\int x \,dx$,
but $\int s \,ds$, with value $\frac 12 s^2$.
Putting the latter in Stieltjes terms, $\int s \,ds$ is the Stieltjes integral $\int_0^t x(s) dx(s)$ where the sample path path or function $x$ is the
identity function $x(s)=s$, $0\leq s\leq t$.

Clearly a Stieltjes integral involving a ``typical'' Brownian path $x(s)$ (which though continuous is, typically, nowhere differentiable)
is a very different beast from a Stieltjes integral involving the straight line path $x(s)=s$.
So in reality it is not surprising that there is a very big difference between the two integrals
\begin{equation}\label{weak stoch int 1}
\int X(s)dX(s)^2, =\frac 12 X(t)^2 - \frac 12 t,\;\;\;\mbox{ and }\;\;\;\int s\,ds, = \frac 12 s^2.
\end{equation}
The first integral typically involves Stieltjes integrals using very complicated and difficult Brownian paths $x(s)$.
It should be distinguished sharply from the more familiar and simpler Stieltjes integrals in which, for instance, the point function
component of the integrand is a continuous function, and the integrator or interval function is formed from increments of a monotone increasing
or bounded variation function.

It is easy to overlook this distinction.
Example 60 of \cite{MTRV} illustrates the potential pitfall. In this Example, $X_T$ is an arbitrary stochastic process
and, with a fixed partition of $T=\,]0,t]$,
$0=\tau_0<\tau_1< \cdots < \tau_m = t$, the function $\sigma(s)$ is constant for $\tau_{j-1} < s \leq \tau_j$. Example 60
claims, in effect, that the stochastic integral $\int_{\tau_{j-1}}^{\tau_j} \sigma(s)dX_s$ exists for each $j$ in the same way that, for constant $\beta$,
$\int_{\tau_{j-1}}^{\tau_j} \beta \;dX_s$ exists and equals $\beta(X(\tau_j) - X(\tau_{j-1}))$.

But Example \ref{Dirichlet 1} above shows that this claim is false. As a step function, $\sigma(\tau_{j-1})$ is not generally equal to the constant
$\beta = \sigma (s)$ when $s>\tau_{j-1}$. So if the sample path $x(s)$ is the Dirichlet function $D(s)$, the Stieltjes integral
$\int_{\tau_{j-1}}^{\tau_j} \sigma(s)dx(s)$ does not exist, and the claim in Example 60 is invalid.

However, if $X_T$ is a Brownian motion process, then each of the significant sample paths $x(s)$ satisfies a condition of uniform continuity.
In that case Example 60 is valid. But it requires some proof, similar to the proof of Theorem 229 on the succeeding page.

So what is truly surprising in (\ref{weak stoch int 1}) is, not that the two integrals give very different results,
but that any convergence at all can be found for the first integral.

Why is this so? This essay has avoided giving any precise meaning to expr\-essions such as $\int_0^t X dX$---or even to
a random variable $X_s$.
But the meaning of the random variable $\int_0^t X dX$ is somehow representative of a Stieltjes-type integral which can be formulated for \textbf{every}
sample path $\{x(s): 0 < s \leq t\}$. These sample paths  may consist of joined-up straight line segments (as in the archetypical
jagged-line Brownian motion diagram), or smooth paths, or everywhere discontinuous paths (like the Dirichlet function). Thus any claim that
all of the separate and individual
Stieltjes integrals $\int_0^t x(s) dx(s)$ of the class of such sample paths $x$---a very large class indeed---have integral values
$\frac 12 x(t)^2 - \frac 12 t$ must be somehow
challenging and dubious.

The integrals $\int_0^t dX(s) = X(t)$, $\int_0^t dx(s) = x(t)$, show that each member of a large class of Stieltjes integrals \textbf{can}
indeed yield a common, single,
simple result. Our discussion of the Riemann sum calculation of these integrals illustrates how this happens:
regardless of the values of $x(s)$
for $s<t$, adding up increments ensures that \textbf{all} values $x(s)$ cancel out, except the terminal value $x(t)$.

Thus, if $f(s)$ takes constant value $\beta$ for $0 \leq s \leq t$, then, for every sample path $x(s)$, the Riemann-Stieltjes
(and Stieltjes-complete) integral $\int_0^t f(s) dx(s)$ exists, and $\int_0^t f(s) dx(s)= \beta x(t)$ (or $\beta(x(t) - x(0))$ if $x(0)\neq 0$.
This is the basis of the claim that  the stochastic integral $\int_0^t f(s) dX(s)$ exists, and is the random variable $\beta X(t)$.

However, Example \ref{Dirichlet 1}
demonstrates that caution must be exercised in pursuing further the logic of Riemann sum
cancelation.
Because, if
the sample path
$x(s), = D(s)$, the expression $f(s)D(I)$ is not integrable on $[0,t]$,
in either the Riemann-Stieltjes sense or the Stieltjes-complete sense,
 even when $f(s)$ is a step function (non-constant).

It is indeed possible to take the Riemann sum
cancelation idea further.  Theorem 229 of \cite{MTRV} shows how this can be done.

But many important stochastic integrands are not actually integrable in the basic sense of the Definition \ref{def Burkill integral}.
If various sample paths $x(s)$ are experimented with in the integral $\int_0^t dX_s^2$, many different results will be found.
So what is the meaning of the result $\int_0^t dX_s^2 = t$?

While, for different sample paths $x$, $\int_0^t dx_s^2$ is not generally convergent to any definite value,
there is a weak sense of convergence of the integral which makes
``$\int_0^t dX_s^2 = t$'' meaningful. Most importantly in this case, the weak limit $t$ is a fixed quantity
rather than a random or unpredictable quantity
such as $x(t)$.
But this question goes beyond the scope of the present essay, whose aim is to explore some of the basic concepts of this subject, and hopefully to
illuminate them a little. A more extensive exploration is given in \cite{MTRV}.


\begin{thebibliography}{99}
\bibitem{Henstock 1955} Henstock, R., \emph{The efficiency of convergence factors for functions of a continuous real variable},
Journal of the London Mathematical Society 30 (1955), 273--286.

\bibitem{MTRV} Muldowney, P., \emph{A Modern Theory of Random Variation, with Applications in Stochastic Calculus,
Financial Mathematics, and Feynman Integration}, Wiley, New York, 2012.
\bibitem{O} {\O}ksendal, B., \emph{Stochastic Differential Equations},   Springer-Verlag, Berlin, 1985.


\end{thebibliography}
\end{document}